\documentclass[11pt]{article}

\font\smallit=cmti10

\usepackage{amssymb,latexsym,amsmath,epsfig,amsthm,mathrsfs,color} 

\makeatletter

\renewcommand\section{\@startsection {section}{1}{\z@}
{-30pt \@plus -1ex \@minus -.2ex}
{2.3ex \@plus.2ex}
{\normalfont\normalsize\bfseries}}

\renewcommand\subsection{\@startsection{subsection}{2}{\z@}
{-3.25ex\@plus -1ex \@minus -.2ex}
{1.5ex \@plus .2ex}
{\normalfont\normalsize\bfseries}}

\renewcommand{\@seccntformat}[1]{\csname the#1\endcsname. }

\makeatother

\begin{document}

\begin{center}
\uppercase{\bf APPEARANCE OF BALANCING AND RELATED NUMBER SEQUENCES IN STEADY STATE PROBABILITIES OF SOME MARKOV CHAINS}
\vskip 20pt
{\bf Asim Patra\footnote{any footnote here}}\\
{\smallit Department of Mathematics, National Institute of technology, Rourkela-769008, Odisha, India}\\
{\tt 515ma3016@nitrkl.ac.in}\\
\vskip 20pt
{\bf Gopal Krishna Panda\footnote{any footnote here}}\\
{\smallit Department of Mathematics, National Institute of technology, Rourkela-769008, Odisha, India}\\
{\tt gkpanda\_nit@rediffmail.com}
\end{center}
\vskip 30pt
\vskip 30pt

2010 Mathematics Subject Classification (A. M. S.):	11B39, 60J05
\vskip 10pt

\centerline{\bf Abstract}
\noindent
Balancing and Lucas-balancing numbers are solutions of a Diophantine equation and satisfy a second order homogeneous recurrence relation. Interestingly, these numbers can be seen as numerators and denominators in the steady state probabilities of a class of transition probability matrices of Markov chains.

\section{Introduction}
A discrete time Markov chain is a stochastic process $ \lbrace X_{k} \rbrace,$ where $k$ runs over non negative integers, such that $Pr\lbrace X_{k+1}=j|X_{k}=i, X_{k-1}=l,..., X_{0}=r\rbrace $ =  $Pr\lbrace X_{k+1}=j|X_{k}=i\rbrace.$ In other words, the Markovian property asserts that, the probability of future behavior of a process depends on the present state of the process and not on the past. The probability $Pr\lbrace X_{k+1}=j|X_{k}=i\rbrace $ not only depends on the initial and final states $i$ and $j$, but it also depends on the time of transition $k.$ When this probability depends only on $i$ and $j$ and not on $k,$ the Markov chain $\lbrace X_{k}\rbrace$ is said to have stationary transition probabilities and in this case  $P_{ij}$=$Pr\lbrace X_{k+1}=j|X_{k}=i\rbrace$ is nothing but the probability of passing from state $i$ to state $j$ in one transition and the matrix  $P=(P_{ij})$ is known as the transition probability matrix (\cite{stochas},\cite{medhi}).
\\
The probability $P_{ij}^{(n)}=Pr\lbrace X_{k+n}=j|X_{k}=i\rbrace$ is called a n-step transition probability. This is the probability of passing from state $i$ to state $j$ in $n$ transitions and the number $P_{ij}^{(n)}$ is the $ij^{th}$ entry of the matrix $\textbf{P}^{n}$(the $n^{th}$ power of the matrix $\textbf{P}$). As $n\rightarrow \infty$, $P_{ij}^{(n)}$ becomes independent of the starting state $i$ and for $j=0,1,2,...$ the numbers $\pi_{j}=\lim_{n\rightarrow \infty} P_{ij}^{(n)}$ are known as the steady state probabilities of the Markov chain $X_{k}.$ The vector $\vec{\pi}=(\pi_{0},\pi_{1},\pi_{2},\cdots)$ satisfies the relationships $\vec{\pi}=\vec{\pi}\textbf{P}$ and $\pi_{0}+\pi_{1}+\pi_{2}+\cdots=1$.
\\
\\
A balancing number, on the other hand, is a natural number $B$ such that $8B^{2}+1 $ is a perfect square and $C=\sqrt{8B^{2}+1} $ is called a Lucas balancing number. The balancing as well as the Lucas-balancing numbers satisfies the same recurrence relation $ x_{n+1}=6x_{n}-x_{n-1}$ with initial terms $ x_{0}=0, x_{1}=1$ and   $x_{0}=1, x_{1}=3$ respectively(see \cite{sqrt},\cite{Thes}). Panda and Rout \cite{recur} studied a class of binary recurrences $x_{n+1}=Ax_{n}-x_{n-1}$  with initial terms $x_{0}=0, x_{1}=1.$ The sequence $\lbrace x_{n}\rbrace_{n=1}^{\infty}$ coincides with the natural numbers when $A=2$ and represents the class of even indexed Fibonacci numbers  when $A=3$; further when $A=6,$ it represents the balancing numbers.

A  sequence closely related to the sequence of balancing numbers is the sequence of cobalancing numbers. A non negative integer $b$ is a cobalancing number if $8b^{2}+8b+1$ is a perfect square and $\sqrt{8b^{2}+8b+1}$ is called a Lucas-cobalancing number. The cobalancing number satisfies the non homogeneous recurrence relation $b_{n+1}=6b_{n}-b_{n-1}+2$ with initial terms $b_{0}=0, b_{1}=0$ where as the Lucas-cobalancing sequence $\lbrace c_{n}\rbrace_{n=1}^{\infty}$ satisfies the same recurrence relation as that of balancing numbers with the initial terms $c_{0}=1, c_{1}=7.$ A very important relationship among balancing and cobalancing number is that the sum of first $n-1$ balancing numbers is equal to half of the $n^{th}$ cobalancing number, that is $\sum_{i=1}^{n-1} B_{i}=b_{n}/2 $(see \cite{recur}).

The Pell sequence is defined by the recurrence relation $P_{n+1}=2P_{n}+P_{n-1},$ $P_{0}=0,P_{1}=1.$ The Pell sequence is related to the balancing sequence by several means. For example, $P_{2n}=2B_{n}, P_{2n+1}=B_{n+1}-B_{n}, n=1,2,\cdots$ For several other relationships among Pell and Balancing sequence the readers are advised to refer to \cite{links}.
It is well known that if $n$ is large, then the $(n+1)^{st}$ Fibonacci number $F_{n+1}$ is approximately $\dfrac{1+\sqrt{5}}{2}$ times the $n^{th}$ Fibonacci numbers $F_{n}$ and the number is known as the golden ratio \cite{golden}. So far as the balancing number is concerned, the $(n+1)^{st}$ balancing number  $B_{n+1}$ is approximately $3+\sqrt{2}$ times the $n^{th}$ balancing numbers $B_{n}[4]$ and the approximation is very sharp. Interestingly $3+2\sqrt{2}=(1+\sqrt{2})^{2}$ and $1+\sqrt{2}$ is known as the silver ratio \cite{silver} and it is interesting to note that for large $n,$ $P_{n+1}$ is approximately equal to $(1+\sqrt{2})P_{n}.$ The approximations hold good even for moderate large value of $n$.

In \cite{matrices}, using matrix
\\
\begin{center}
$B=\begin{bmatrix}
6 &  1\\ 
-1 & 0
\end{bmatrix}=\begin{bmatrix}
B_{2} &  B_{1}\\ 
-B_{1} & B_{0}
\end{bmatrix}
,$
\end{center}

the balancing numbers have been calculated using the relation
\begin{center}
$B^{n}=\begin{bmatrix}
B_{n+1} &  B_{n}\\ 
-B_{n} & B_{n-1}
\end{bmatrix}.$
\end{center}
Hlynka and Sajobi\cite{markv} established the presence of Fibonacci numbers in numerators and denominators of the stesdy state probabilities of a particular class of Markov chains.Motivated by their work, we construct a class of Markov chains such that their steady state distributions involve balancing, Lucas-balancing and balancing-like numbers.

\section{BALANCING NUMBERS IN STEADY STATE PROBABILITIES}
Consider a village which can accommodate not more than a population of size $n.$ Let $X_{k}$ be the population of the village to increase by one for each birth, decrease by one for each death and become zero when the whole population is migrated to a different destination. Then $\lbrace X_{k}\rbrace_{k=0}^{\infty}$ can be viewed as a Markov Chain and let us assign the transition probabilities $P_{i,i+1}=P_{i,i-1}=\dfrac{1}{6}$ if $1\leq i\leq n-2, P_{01}=\dfrac{1}{6}, P_{00}=P_{10}=P_{n-1,0}=\dfrac{5}{6}, P_{i0}=\dfrac{2}{3}$ if $2\leq i\leq n-2$ and $P_{ij}=0$ otherwise.
\\ Thus the transition probability  matrix is given by

\begin{center}
$\textbf{P}$= 
$\begin{bmatrix}
&5/6&  1/6&  0&   0&  0&  \cdots&  0& 0&\\ 
&5/6&  0&  1/6& 0&  0&  \cdots&  0&  0&  \\ 
&2/3&  1/6& 0&  1/6&  0&  \cdots&  0&  0&\\
&2/3&  0& 1/6&  0&  1/6&  \ddots&  0&  0&\\

&\vdots&  \vdots& \ddots&  \ddots&  \ddots&  \ddots&  \ddots&  \ddots&\\
&2/3&  0& 0&  0&  0&  \ddots&  0&  1/6&\\
&5/6&  0& 0&  0&  0&  \ddots&  1/6&    0&\\
\end{bmatrix}$ ~~~~~~~~~~                (1)
\end{center}

and the vector $ \vec{\pi}=(\pi_{0},\pi_{1}, \pi_{2},\cdots,\pi_{n-1})$ representing the steady state probabilities of $\lbrace X_{k}\rbrace_{k=0}^{\infty}$ can be calculated from the relation 
$\vec{\pi}=\vec{\pi}\textbf{P}$ subject to the condition $\pi_{0}+\pi_{1}+\cdots+\pi_{n-1}=1.$
The following theorem shows that the steady state probabilities are functions of balancing and cobalancing numbers.\\

\textbf{Theorem 2.1.}~~The steady state probability  vector corresponding to the transition probability matrix $\textbf{P}$ given in (1) is $\vec{\pi}= \left(\frac{2B_{n}}{b_{n+1}},\frac{2B_{n-1}}{b_{n+1}},\cdots,\frac{2B_{1}}{b_{n+1}}\right).$
\\

\textbf{Proof}. The steady state probability  vector corresponding to the transition probability matrix $\textbf{P}$ can be obtained from the relationship $\vec{\pi}=\vec{\pi}\textbf{P}$ leading to the system of equations(in reverse order)
\begin{center}
$\pi_{n-1}=\dfrac{1}{6}\pi_{n-2},$ $\pi_{i}=\dfrac{1}{6}\pi_{i-1}+\dfrac{1}{6}\pi_{i+1}, i=1,2,\cdots, n-2.$
\end{center}
These equations are known as balance equations. On rearranging them, we get
\begin{center}
$\pi_{n-2}=6\pi_{n-1},$ $\pi_{i-1}=6\pi_{i}-\pi_{i+1}, i=1,2,\cdots, n-2.$~~~~~~~~~~~~~~~~(2)
\end{center}

We ignore the equation $\pi_{0}=\dfrac{5}{6}\pi_{0}+\dfrac{5}{6}\pi_{1}+\dfrac{2}{3}(\pi_{2}+\cdots+\pi_{n-2}+\dfrac{5}{6}\pi_{n-1}$ since we have the additional equation  $\pi_{0}+\pi_{1}+\cdots+\pi_{n-1}=1.$ Setting $\pi_{n-1}=k,$ we can rewrite the system of equations in (2) as 
\begin{center}
$\pi_{n-1}=k=kB_{1},$ $\pi_{n-2}=6k=kB_{2},$ $\pi_{i-1}=6\pi_{i}-\pi_{i+1}, i=1,2,\cdots, n-2.$ 
\end{center}
We will show that $\pi_{i}=kB_{n-i}$ for $i=0,1,2,\cdots, n-1$ using mathematical induction.
\\
~~~~~Since $\pi_{n-1}=k=kB_{1},$ $\pi_{n-2}=6k=kB_{2}$, the assertion is true for $i=0,1.$ Assume that the assertion is true for $i=j\leq n-1.$ Since the assertion is assumed to be true for $i=j,$ $\pi_{j}=kB_{n-j}$. In  view of the equation $\pi_{j-1}=6\pi_{j}+\pi_{j+1},$ we have 
\begin{center}
$\pi_{j+1}=6kB_{n-j}-kB_{n-j+1}=k(6B_{n-j}-B_{n-j+1})$ 
\end{center}
and since $B_{j+1}=6B_{j}-B_{j-1},$it follows that $\pi_{j+1}=kB_{n-j-1}=kB_{n-(j+1)}$ and the assertion is true for $i=j+1\leq n-1.$ Further, $\pi_{0}+\pi_{1}+\cdots+\pi_{n-1}=1$ implies that $k=\frac{1}{\sum_{l=1}^{n}{B_{l}}}.$ Hence,
\begin{center}
 $\pi_{i}=\frac{2B_{n-i}}{\sum_{l=1}^{n}{B_{l}}}.$ for $i=0,1,\cdots,n-1.$     ~~~~~~~~~~~~~~~~~~~~~~~(3)
\end{center}
Since $\sum_{l=1}^{n}{B_{l}}=\frac{b_{n+1}}{2},$ the proof is complete.~~$\blacksquare$

\textbf{Remark:} Since $2B_{n}=P_{2n}$, (3) can be rewritten as $\pi_{i}=\frac{2P_{2n}}{\sum_{l=1}^{n}{B_{l}}}$ for $i=0,1,\cdots,n-1.$

The $n\times n$ transition probability matrix $\textbf{P}$ defined in (1) results in the steady state probabilities (3). In the following theorem, we will show that a class of transition probability matrices in which $\textbf{P}$ defined in (1) is a member, result in the steady state probabilities (3).
\\

\textbf{Theorem 2.2.}~~~If $0<q\leq 1/6 ,$ the transition probability matrix
\begin{center}
$\textbf{P}(q)$= 
$\begin{bmatrix}
&1-q&  q&  0&   0&  0&  \cdots&  0& 0&\\ 
&5q&  1-6q&  q& 0&  0&  \cdots&  0&  0&  \\ 
&4q&  q& 1-6q&  q&  0&  \cdots&  0&  0&\\
&4q&  0& q&  1-6q&  q&  \ddots&  0&  0&\\

&\vdots&  \vdots& \ddots&  \ddots&  \ddots&  \ddots&  \ddots&  \ddots&\\
&4q&  0& 0&  0&  0&  \ddots&  1-6q&  q&\\
&5q&  0& 0&  0&  0&  \ddots&  q&    1-6q&\\
\end{bmatrix}$~~~~~~~~~~~~~~~~~~~~~(4)
\end{center}
results in the same steady state probabilities (3).

\textbf{Proof}. Let $q$ be any number with $0<q\leq 1/6 .$ The steady state probability vector corresponding to the transition probability matrix $\textbf{P}(q)$ can be obtained from $\vec{\pi}=\vec{\pi}\textbf{P}. $ This leads to the system of equations(in reverse order) 
\begin{center}
$\pi_{n-1}=q\pi_{n-2}+(1-6q)\pi_{n-1},$ $\pi_{i}=q\pi_{i-1}+(1-6q)\pi_{i}+q\pi_{i+1}, i=1,2,\cdots, n-2.$~~~~~~~~~~~~~~(5)
\end{center}
As usual, the equation $\pi_{0}=(1-q)\pi_{0}+5q\pi_{1}+4q(\pi_{2}+\cdots+\pi_{n-2})+5q\pi_{n-1}$  is ignored because of the additional equation $\pi_{0}+\pi_{1}+\cdots+\pi_{n-1}=1$ satisfied by the steady state probabilities. On simplifying (5), we get

\begin{center}
$\pi_{n-2}=6\pi_{n-1},$  $\pi_{i-1}=6\pi_{i}-\pi_{i+1}, i=1,2,\cdots, n-2.~~~~~~~~~~~~~~~~~~(6)$
\end{center}
Since the system of equations in (6) is essentially same as those in (2), the conclusion of the theorem follows.~~$\blacksquare$

We noticed in the remark following the proof of Theorem 2.1 that the transition probability matrix $\textbf{P}$  given in (1) has steady state probabilities that are rational functions of even indexed Pell numbers. We next present a transition probability matrix which is not much different from $\textbf{P}$ given in (1), but its steady state probabilities can be expressed as rational function of consecutive Pell numbers.\\

\textbf{Theorem 2.3.} The steady state probability vector of the transition probability matrix
\begin{center}
$\textbf{P}$= 
$\begin{bmatrix}
&5/6&  1/6&  0&   0&  0&  \cdots&  0& 0&\\ 
&5/6&  0&  1/6& 0&  0&  \cdots&  0&  0&  \\ 
&2/3&  1/6& 0&  1/6&  0&  \cdots&  0&  0&\\
&2/3&  0& 1/6&  0&  1/6&  \ddots&  0&  0&\\

&\vdots&  \vdots& \ddots&  \ddots&  \ddots&  \ddots&  \ddots&  \ddots&\\
&2/3&  0& 0&  0&  0&  \ddots&  0&  1/6&\\
&2/3&  0& 0&  0&  0&  \ddots&  1/6&    1/6&\\
\end{bmatrix}$~~~~~~~~~~~~~~~~~~~~~~~~~(7)
\end{center}
results in the same steady state probabilities (3).

$(P_{i,i+1}=P_{i,i-1}=\dfrac{1}{6}$ if $1\leq i\leq n-1, P_{01}=\dfrac{1}{6}, P_{00}=P_{10}=\dfrac{5}{6}, P_{i0}=\dfrac{2}{3}$ if $2\leq i\leq n-1$ and $P_{ij}=0$ otherwise) is $\vec{\pi}=(\pi_{0},\pi_{1},\pi_{2},\cdots,\pi_{n-1})$ where $\pi_{i}=\dfrac{B_{n+i}-B_{n+i-1}}{B_{n}}$ for $i=0,1,\cdots,n-1$.

\textbf{Proof}. The balance equation for the steady state probabilities corresponding to the transition probability matrix \textbf{P} are given by 
\begin{center}
$\pi_{n-1}=\dfrac{1}{6}\pi_{n-2}+\dfrac{1}{6}\pi_{n-1},$ $\pi_{i}=\dfrac{1}{6}\pi_{i-1}+\dfrac{1}{6}\pi_{i+1}, i=1,2,\cdots, n-2.$
\end{center}
 On rearrangement, we get
\begin{center}
$\pi_{n-2}=5\pi_{n-1},$ $\pi_{i-1}=6\pi_{i}-\pi_{i+1}, i=1,2,\cdots, n-2.$
\end{center}
Letting $\pi_{n-1}=k,$ we see that $\pi_{n-2}=5k.$ Observe that $\pi_{n-i}=k\left(B_{i}-B_{i-1}\right)$ is true for $i=1,2.$ Assuming that the assertion is true for $i=r,$ then the balance equation 
$\pi_{n-r-1}=6\pi_{n-r}-pi_{n-r+1}$ implies that
\begin{center}
$\pi_{n-r-1}=6k(6B_{r}-B_{r-1})-k(B_{r-1}-B_{r-2})=k(B_{r+1}-B_{r})$ 
\end{center}
and hence the assertion is true for $i=r+1.$ Since $\pi_{0}+\pi_{1}+\cdots+\pi_{n-1}=1,$ the conclusion of the theorem follows.~~$\blacksquare$

\textbf{Remark}: Since $P_{2n}=2B_{n}$ and $P_{2n-1}=B_{n}-B_{n-1},$ we can express the steady state probabilities in the previous theorem as $\pi_{i}=\dfrac{2P_{2n-1}}{P_{2n}}, i=0,1,\cdots,n-1.$

\section{LUCAS BALANCING NUMBERS IN STEADY STATE PROBABILITIES}
In the previous section we established the appearance of balancing numbers in the steady state probabilities of some Markov chains. In the present section, we introduce some finite state Markov chains whose steady state probability distributions are rational functions of Lucas-balancing numbers.
\\
To start with, let $\lbrace X_{k}\rbrace_{k=0}^{\infty}$ be a Markov chain having the state space $\lbrace 0,1,2,\cdots,n-1\rbrace$ and transition probability matrix

\begin{center}
$\textbf{P}$= 
$\begin{bmatrix}
&5/6&  1/6&  0&   0&  0&  \cdots&  0& 0&\\ 
&5/6&  0&  1/6& 0&  0&  \cdots&  0&  0&  \\ 
&2/3&  1/6& 0&  1/6&  0&  \cdots&  0&  0&\\
&2/3&  0& 1/6&  0&  1/6&  \ddots&  0&  0&\\

&\vdots&  \vdots& \ddots&  \ddots&  \ddots&  \ddots&  \ddots&  \ddots&\\
&2/3&  0& 0&  0&  0&  \ddots&  0&  1/6&\\
&1/3&  0& 0&  0&  0&  \ddots&  1/6&    1/2&\\
\end{bmatrix}$~~~~~~~~~~~~~~~~~~~~(8)
\end{center}

More specifically  $\textbf{P}=(\textbf{P}_{ij}),$  $P_{i,i+1}=P_{i,i-1}=\dfrac{1}{6}$ if $1\leq i\leq n-2, P_{01}=\dfrac{1}{6}, P_{00}=P_{10}=\dfrac{5}{6}, P_{i0}=\dfrac{2}{3}$ if $2\leq i\leq n-2,$ $P_{n-1,0}=\dfrac{1}{3}, P_{n-1,n-1}=\dfrac{1}{2}$ and $P_{ij}=0$ otherwise. In the following theorem, we will see that the steady state probability vector corresponding to $\textbf{P}$ used in Theorem 2.1; the only difference is that the balancing numbers appearing in the steady state probabilities in Theorem 2.1 will be replaced by Lucas-balancing numbers.\\

\textbf{Theorem 3.1.} The steady state probability vector corresponding to the transition probability matrix \textbf{P} given in (8) is  $\vec{\pi}=(\pi_{0},\pi_{1},\pi_{2},\cdots,\pi_{n-1})$ where $\pi_{i}=\frac{C_{n-i}}{\sum_{i=1}^{n}{C_{i}}}$ for $i=0,1,\cdots,n-1$ and $C_{i}$ denotes the $i^{th}$ Lucas-balancing numbers.\\

\textbf{Proof}.    The balance equation relating the steady state probabilities associated with $\textbf{P}$ are given by 
\begin{center}
$\pi_{n-1}=\dfrac{1}{6}\pi_{n-2}+\dfrac{1}{2}\pi_{n-1},$ $\pi_{i}=\dfrac{1}{6}\pi_{i-1}+\dfrac{1}{6}\pi_{i+1}, i=1,2,\cdots, n-2.$
\end{center}
 On rearrangement, we get
 \begin{center}
 $\pi_{n-2}=3\pi_{n-1},$ $\pi_{i-1}=6\pi_{i}-\pi_{i+1}, i=1,2,\cdots, n-2.$\end{center}
Letting $\pi_{n-1}=k,$ we see that $\pi_{n-2}=3k.$ Observe that $\pi_{n-i}=kC_{n-i}$ is true for $i=1,2.$ Assuming that the assertion is true for $i=r,$ then the balance equation $\pi_{n-r-1}=6\pi_{n-r}-\pi_{n-r+1}$ implies that
\begin{center}
$\pi_{n-r-1}=6kC_{r}-kC_{r-1}=kC_{r+1}$ 
\end{center}
and hence the assertion is true for $i=r+1.$ Since $\pi_{0}+\pi_{1}+\cdots+\pi_{n-1}=1,$ the conclusion of the theorem follows.~~$\blacksquare$
\\

The $n\times n$ transition probability matrix $\textbf{P}$ defined in (8) results in the steady state probabilities  $\pi_{i}=\frac{C_{n-i}}{\sum_{i=1}^{n}{C_{i}}}$ for $i=0,1,\cdots,n-1.$ In the following theorem, we will show that a class of transition probability matrices in which $\textbf{P}$ defined in (8) is a member, also result in  the same steady state probabilities.\\

\textbf{Theorem 3.2.}~~If $q$ is a real number such that $0<q\leq 1/6 ,$ then the transition probability matrix
\begin{center}
$\textbf{P(q)}$= 
$\begin{bmatrix}
&1-q&  q&  0&   0&  0&  \cdots&  0& 0&\\ 
&5q&  1-6q&  q& 0&  0&  \cdots&  0&  0&  \\ 
&4q&  q& 1-6q&  q&  0&  \cdots&  0&  0&\\
&4q&  0& q&  1-6q&  q&  \ddots&  0&  0&\\

&\vdots&  \vdots& \ddots&  \ddots&  \ddots&  \ddots&  \ddots&  \ddots&\\
&4q&  0& 0&  0&  0&  \ddots&  1-6q&  q&\\
&2q&  0& 0&  0&  0&  \ddots&  q&    1-3q&\\
\end{bmatrix}$  ~~~~~~~~~~~~~~~~~~(9)
\end{center}
results in the same steady state probabilities obtained in Theorem 3.1.

The proof of the above theorem is similar to that of Theorem 2.2 and hence it is omitted.
\\\\
Lucas-cobalancing numbers are associated with cobalancing numbers in the same manner the Lucas-balancing numbers are associated with balancing numbers. However, as we have discussed in Section 1, their recurrence relation is identical to that of balancing numbers.
In the following theorem, we consider a finite state Markov chain whose steady state probabilities are rational functions of Lucas-cobalancing numbers. Since the proof is similar to that of Theorem 3.1, we prefer to omit it.\\

\textbf{Theorem 3.3.}~~The steady state probability vector corresponding to the transition probability matrix 
\begin{center}
$\textbf{P}$= 
$\begin{bmatrix}
&5/6&  1/6&  0&   0&  0&  \cdots&  0& 0&\\ 
&5/6&  0&  1/6& 0&  0&  \cdots&  0&  0&  \\ 
&2/3&  1/6& 0&  1/6&  0&  \cdots&  0&  0&\\
&2/3&  0& 1/6&  0&  1/6&  \ddots&  0&  0&\\

&\vdots&  \vdots& \ddots&  \ddots&  \ddots&  \ddots&  \ddots&  \cdots&\\
&2/3&  0& \ddots&  \ddots&  \ddots&  \ddots&  \ddots&  \cdots&\\
&16/21&  0& 0&  0&  0&  \ddots&  0&  1/14&\\
&1/3&  0& 0&  0&  0&  \ddots&  1/6&    1/2&\\
\end{bmatrix}$~~~~~~~~~~~~~~~~~~~(10)
\end{center}
  $(P_{i,i+1}=P_{i,i-1}=\dfrac{1}{6}$ if $1\leq i\leq n-2, P_{01}=\dfrac{1}{6}, P_{00}=P_{10}=\dfrac{5}{6}, P_{i0}=\dfrac{2}{3}$ if $2\leq i\leq n-3,$ $P_{n-2,0}=\dfrac{16}{21},$  $P_{n-1,0}=\dfrac{1}{3},$  $P_{n-2,n-1}=\dfrac{1}{14}$   $P_{n-1,n-1}=\dfrac{1}{2}$ and $P_{ij}=0$ otherwise), is  $\vec{\pi}=(\pi_{0},\pi_{1},\pi_{2},\cdots,\pi_{n-1})$ where $\pi_{i}=\dfrac{c_{n-i}}{\sum_{i=1}^{n}{c_{i}}}$ for $i=0,1,\cdots,n-1.$  or $\vec{\pi}=\left( \dfrac{c_{n}}{\sum_{i=1}^{n}{c_{i}}}, \dfrac{c_{n-1}}{\sum_{i=1}^{n}{c_{i}}} , \cdots, \dfrac{c_{1}}{\sum_{i=1}^{n}{c_{i}}} \right)$ where $c_{n}$ denotes the $n^{th}$ Lucas-cobalancing number.

\section{SILVER RATIO IN STEADY STATE PROBABILITIES OF MARKOV CHAINS WITH INFINITE STATE SPACES}

In the last two sections, we have studied steady state probabilities of Markov chains with finite state spaces. In this section, we consider a Markov chain having the infinite state space $\lbrace 0,1,2,\cdots \rbrace$ and the transition probability matrix

\begin{center}
$\textbf{P}$= 
$\begin{bmatrix}
&5/6&  1/6&  0&   0&  0&  \cdots&  \\ 
&5/6&  0&  1/6& 0&  0&  \cdots&   \\ 
&2/3&  1/6& 0&  1/6&  0&  \cdots& \\
&2/3&  0& 1/6&  0&  1/6&  \ddots& \\

&\vdots&  \vdots& \ddots&  \ddots&  \ddots&  \ddots& \\
\end{bmatrix}.$~~~~~~~~~~~~~(11)
\end{center}

More specifically, \textbf{P}= $\textbf{P}_{ij}$  with $P_{i,i+1}=\dfrac{1}{6}$ if $i\geq 1, P_{10}=\dfrac{5}{6}, P_{i0}=\dfrac{2}{3}$ if $i\geq 2,$ and $P_{ij}=0$ otherwise. In the following theorem we find the steady state probabilities corresponding to $\textbf{P}$  in terms of odd indexed Pell numbers. To do so, we will study the limiting behavior of steady state probabilities obtained in Section 2.\\

\textbf{Theorem 4.1.} ~~~The steady state probability vector corresponding to the transition probability matrix $\textbf{P}$ given in (11) is $\vec{\pi}=(\pi_{0},\pi_{1},\pi_{2},\cdots)$ where $\pi_{i}=\beta^{i}-\beta^{i+1}, i=0,1,2,\cdots$ and $\beta=3-2\sqrt{2}.$
\\
\textbf{Proof.} ~~Using the identity $\vec{\pi}=\vec{\pi} \textbf{P},$  the balance equations for the calculation of steady state probabilities are given by 

\begin{center}
$\pi_{0}=\dfrac{5}{6}\pi_{0}+\dfrac{5}{6}\pi_{1}+\dfrac{2}{3}(\pi_{2}+\pi_{3}+\cdots),$ $\pi_{i}=\dfrac{1}{6} \pi_{i-1}+\dfrac{1}{6} \pi_{i+1}$ for $i \geq 1.$
\end{center}
On simplification, we get $\pi_{1}=5\pi_{0}-4, \pi_{2}=29\pi_{0}-24 , \pi_{3}=169\pi_{0}-140 $ and using mathematical induction, one can see that
\begin{center}
$\pi_{i}=(B_{i+1}-B_{i})\pi_{0}-4B_{i}, i=1,2,\cdots$~~~~~~~~~~~~~(12)
\end{center}
Solving the infinite system(12) is not easy. Instead, we observe that transition probability matrix \textbf{P} given in (11) is a limiting case of the $n \times n$ transition probability matrix \textbf{P} that appear in Theorem 2.3. In the proof of Theorem 2.3, we have noticed that $\pi_{0}=\frac{B_{n}-B_{n-1}}{B_{n}}.$ Hence, in this case, we have
\begin{center}
$\pi_{0}=lim_{n\rightarrow \infty} \dfrac{B_{n}-B_{n-1}}{B_{n}}$=$1-lim_{n\rightarrow \infty} \dfrac{B_{n-1}}{B_{n}}=1-(3-2\sqrt{2})=1-\beta$ \end{center}

where $\beta=3-2\sqrt{2}.$ Similarly,
\begin{center}
$\pi_{1}=lim_{n\rightarrow \infty} \dfrac{B_{n-1}-B_{n-2}}{B_{n}}$=$lim_{n\rightarrow \infty} \left(\dfrac{B_{n-1}}{B_{n}}- \dfrac{B_{n-2}}{B_{n-1}}.\dfrac{B_{n-1}}{B_{n}}\right)=\beta-\beta^{2},$

$\pi_{2}=lim_{n\rightarrow \infty} \dfrac{B_{n-2}-B_{n-3}}{B_{n}}$=$\beta^{2}-\beta^{3},$
\end{center}

and in general
\begin{center}
$\pi_{i}=lim_{n\rightarrow \infty} \dfrac{B_{n-i-1}-B_{n-i}}{B_{n}}$=$\beta^{i}-\beta^{i+1},i=0,1,2,\cdots$~~~$\blacksquare$
\end{center}
It is easy to verify that if the transition probability matrix \textbf{P} used in Theorem 4.1 is replaced by 

\begin{center}
$\textbf{P}$= 
$\begin{bmatrix}
&1-q&  q&  0&   0&  0&  \cdots&  \\ 
&5q&  1-6q&  q& 0&  0&  \cdots&   \\ 
&4q&  q& 1-6q&  q&  0&  \cdots& \\
&4q&  0& q&  1-6q&  q&  \cdots& \\

&\vdots&  \vdots& \ddots&  \ddots&  \ddots&  \ddots& \\
\end{bmatrix}$~~~~~~~~~~~~~~(13)
\end{center}

where $0<q\leq 1/6,$ the conclusion of Theorem 4.1 remain unchanged.\\
\textbf{Remark:} Observe that $\beta=\dfrac{1}{(1+\sqrt{2})^{2}}$ and ratio $1+\sqrt{2}:1$ is nothing but the silver ratio discussed in Section 1.\\

~~~~~~~The following corollary, which is a consequence of the above theorem, establishes an interesting relationship among balancing numbers and $\beta.$ Observe that the proof of the identity remains probabilistic.

\textbf{Corollary 4.2.} $\beta^{n+1}=\beta B_{n+1}-B_{n}, n=1,2,\cdots,$ where   $\beta=3-2\sqrt{2}.$

\textbf{Proof.} In view of Theorem 4.1, we have $\pi_{0}=1-\beta$ and  $\beta^{i}-\beta^{i+1}=\beta^{i}(1-\beta).$ Also from (12) we have
\begin{center}
$\pi_{i}=(B_{i+1}-B_{i})\pi_{0}-4B_{i},i=1,2,\cdots$
\end{center}
Substituting the value of $\pi_{0},$ we get 
\begin{center}
$\beta^{i}(1-\beta)=(B_{i+1}-B_{i})(1-\beta)-4B_{i}.$
\end{center}
Thus,
\begin{center}
$\beta^{i}=(B_{i+1}-B_{i})-\dfrac{4B_{i}}{1-\beta}=B_{i+1}-(3+2\sqrt{2})B_{i}=B_{i+1}-\dfrac{B_{i}}{\beta}$
\end{center}
from which the conclusion of the corollary follows.~~$\blacksquare$

\section{BALANCING-LIKE NUMBERS IN THE STEADY STATE PROBABILITIES OF MARKOV CHAINS}

In Section 1, we have seen that if $A>2$ is a fixed positive integer, then the sequence $\lbrace x_{n} \rbrace_{n=1}^{\infty}$ defined recursively by $x_{n+1}=Ax_{n}-x_{n-1}$ with initial terms $x_{0}=0, x_{1}=1$ is called a balancing-like sequence and this sequence serve as a generalization of the balancing sequence. In this section, we consider a finite state Markov chain whose steady state probabilities involve balancing-like numbers.
\\
To start with we consider the following $n \times n$ transition probability matrix
 
\begin{center}
$\textbf{P}(A)$= 
$\begin{bmatrix}
&1-\frac{1}{A}&  \frac{1}{A}&  0&   0&  0&  \cdots&  0& 0&\\ 
&1-\frac{1}{A}&  0&  \frac{1}{A}& 0&  0&  \cdots&  0&  0&  \\ 
&1-\frac{2}{A}&  \frac{1}{A}& 0&  \frac{1}{A}&  0&  \cdots&  0&  0&\\
&1-\frac{2}{A}&  0& \frac{1}{A}&  0&  \frac{1}{A}&  \ddots&  0&  0&\\

&\vdots&  \vdots& \ddots&  \ddots&  \ddots&  \ddots&  \ddots&  \ddots&\\
&1-\frac{2}{A}&  0& 0&  0&  0&  \ddots&  0&  \frac{1}{A}&\\
&1-\frac{1}{A}&  0& 0&  0&  0&  \ddots&  \frac{1}{A}&    0&\\
\end{bmatrix}$~~~~~~~~~~~~~~~~(14)
\end{center}

To be more specific, $P_{i,i+1}=P_{i,i-1}=\frac{1}{A}$ if $1\leq i\leq n-2, P_{01}=\frac{1}{A}, P_{00}=P_{10}=P_{n-1,0}=1-\frac{1}{A}, P_{i0}=1-\frac{2}{A}$ if $2\leq i\leq n-2$ and $P_{ij}=0$ otherwise.

In the following theorem we will see that the steady state probabilities corresponding to $\textbf{P}(A)$ are  rational functions of the first $n$ balancing-like numbers.\\

\textbf{Theorem 5.1.}~~ The steady state probability vector corresponding to the $n \times n$ transition probability matrix $\textbf{P}(A)$ in (14) is $\vec{\pi}=(\pi_{0},\pi_{1},\pi_{2},\cdots,\pi_{n-1})$ where $\pi_{i}=\dfrac{x_{n-i}}{\sum_{i=1}^{n} x_{i}}$ for $i=0,1,\cdots,n-1.$
\\

\textbf{Proof.} The proof of this theorem is similar to that of Theorem 2.1 and hence it is omitted.~~$\blacksquare$
\\
It is important to note that if $q$ is any real number such that $0<q\leq \frac{1}{A},$ then the transition probability matrix

\begin{center}
$\textbf{P}(q)$= 
$\begin{bmatrix}
&1-q&  q&  0&   0&  0&  \cdots&  0& 0&\\ 
&(A-1)q& 1-Aq&  q& 0&  0&  \cdots&  0&  0&  \\ 
&(A-2)q&  q& 1-Aq&  q&  0&  \cdots&  0&  0&\\
&(A-2)q&  0& q&  1-Aq&  q&  \ddots&  0&  0&\\

&\vdots&  \vdots& \ddots&  \ddots&  \ddots&  \ddots&  \ddots&  \ddots&\\
&(A-2)q&  0& 0&  0&  0&  \ddots&  1-Aq&  q&\\
&(A-1)q&  0& 0&  0&  0&  \ddots&  q&    1-Aq&\\
\end{bmatrix}$
\end{center}

gives the steady state probabilities as stated in Theorem 5.1.


\section{CONCLUSION}
In this work, we established the appearance of balancing and related numbers sequence in the steady state probabilities of some Markov chains.  We also noticed that, in many instances, a class of transition probability matrices gives rise to same steady state probabilities. Using the balance equations, we also derived an identity relating the balancing numbers and the silver ratio. Some problems in this area are still open. We encourage the readers to construct transition probability matrices whose steady state vectors would explore some other number sequences.   In this process, they may be able to prove some identities using the balance equations.

\providecommand{\href}[2]{#2}\begingroup
\end{document}